\documentclass{fic-l}       

\usepackage{amssymb}
\usepackage{latexsym}

\newtheorem{thm}{Theorem}[section]
\newtheorem{pro}[thm]{Proposition}
\newtheorem{lem}[thm]{Lemma}

\theoremstyle{definition}
\newtheorem{defi}[thm]{Definition}
\newtheorem{exa}[thm]{Example}

\theoremstyle{remark}

\numberwithin{equation}{section}

\newcommand{\ncm}{\newcommand}

\ncm{\beq}{\begin{equation}}
\ncm{\eeq}{\end{equation}}
\ncm{\bea}{\begin{eqnarray}}
\ncm{\eea}{\end{eqnarray}}
\ncm{\beanon}{\begin{eqnarray*}}
\ncm{\eeanon}{\end{eqnarray*}}

\def\M{\mathsf{M}}
\def\rep{\mathsf{rep}\,}
\def\End{\mbox{\rm End}\,}

\def\Hom{\mbox{\rm Hom}\,}

\def\id{\mbox{\rm id}\,}
\def\Center{\mbox{\rm Center}\,}

\def\B{\mathcal{B}}
\def\C{\mathcal{C}}
\def\D{\mathcal{D}}
\def\o{\otimes}    
\def\oo{\Box\,}    
\def\x{\times}     
\def\bra{\langle}
\def\ket{\rangle}
\def\under{\mbox{\rm\_}\,}
\ncm{\rarr}[1]{\stackrel{#1}{\longrightarrow}}
\ncm{\larr}[1]{\stackrel{#1}{\longleftarrow}}
\def\cop{\Delta}
\def\eps{\varepsilon}
\def\1{_{(1)}}
\def\2{_{(2)}}
\def\3{_{(3)}}
\def\tr{\mbox{\rm tr}}
\def\PL{\pi^{\scriptscriptstyle L}}
\def\PR{\pi^{\scriptscriptstyle R}}

\def\duA{\hat A}
\def\la{\!\rightharpoonup\!}
\def\ra{\!\leftharpoonup\!}
\def\du1{\hat 1}
\def\iso{\rarr{\sim}}
\def\d{{\mathbf d}}
\def\iL{\,\imath_{\scriptscriptstyle L}\,}
\def\iR{\,\imath_{\scriptscriptstyle R}\,}
\def\cL{c_{\scriptscriptstyle L}\,}
\def\cR{c_{\scriptscriptstyle R}\,}

\begin{document}

\title{Finite quantum groupoids and inclusions of finite type} 
\dedicatory{Dedicated to Sergio Doplicher and John E. Roberts on the
occasion of their 60th birthday} 

\author[K. Szlach\'anyi]{Korn\'el Szlach\'anyi}
\address{Research Institute for Particle and Nuclear Physics, Budapest\\    
H-1525 Budapest, P. O. Box 49, Hungary}
\email{szlach@rmki.kfki.hu} 
\thanks{Supported in part by the Hungarian Scientific Research Fund,
OTKA--T 030099} 

\maketitle

In the last decade various motivations coming from low dimensional
quantum field theory, operator algebra, and 
Poisson geometry have lead to the introduction of 
a new notion of symmetry that generalizes both quantum
groups and classical groupoid algebras. The slightly different definitions
\cite{Mal, Yam, Hay, BSz, Lu, Val1, BNSz, Nill, EV} 
all share in the property that a "quantum groupoid" $A$
contains two antiisomorphic canonical
subalgebras: The source subalgebra $A^R$ and the target 
subalgebra $A^L$. They reduce to the scalars in the case of a Hopf algebra
and they are the carrying space of the trivial representation in the monoidal
category $\M_A$ of $A$-modules. In the groupoid interpretation $A^L$
and $A^R$ are non-commutative analogues of the algebra of functions on the
space of units.  The various definitions of quantum
groupoid differ in the size and in commutativity of these subalgebras. The
most general among them is Lu's {\em Hopf algebroid} \cite{Lu}, while the
($C^*$)-{\em weak Hopf algebra} of \cite{BSz, BNSz} captures the most general
"finite quantum groupoid" (see \cite{NV4} for a review) which has the extra
beauty of being selfdual, just like a finite Abelian group or a finite
dimensional Hopf algebra. Yamanouchi's generalized Kac algebra and Hayashi's
face algebra are special cases of the latter. In \cite{EV}, generalizing
\cite{Val1}, Enock and Vallin introduced the notion of a Hopf bimodule which
corresponds to Lu's bialgebroid in the von Neumann algebraic framework. 
In \cite{E2} Enock constructs an antipode for Hopf bimodules making use
of modular theory.

The Doplicher-Roberts duality 
theorem \cite{DR} characterizes the symmetric monoidal Abelian $C^*$-categories
with irreducible monoidal unit as representation categories of (uniquely
determined) compact groups. In this way it provides an intrinsic definition of
internal symmetry of DHR sectors of quantum field theories in spacetime
dimension greater than 2. In dimension 2 no analogue result is known. 
In this respect
the significance of $C^*$-weak Hopf algebras is two-fold.
They have representation categories such that 
\begin{enumerate}
\item[i)] the intrinsic (categorical)
dimensions of the objects (in the sense of \cite{LR}) are not integers
\item[ii)] and the monoidal unit is reducible allowing different 
irreducible "vacuum representations" \cite{BSz2}.
\end{enumerate}
Property (i) offers a possibility that the braided $C^*$-categories found in 
conformal field theory models are equivalent to representation categories of
$C^*$-weak Hopf algebras, although uniqueness cannot be expected.
Even non-braided $C^*$-categories are included, therefore (ii) suggests that
topological soliton sectors can also be described by weak Hopf symmetry.  

The universal problem to which the answer is a unique quantum groupoid is not
known. But inclusions of (unital $C^*$-, von Neumann) algebras $N\subset M$
are very close to that. In \cite{NSzW} the regular action of a $C^*$-weak Hopf
algebra $A$ on a von Neumann algebra $M$ has been defined. This is a kind of
Galois action which allows $M^A\subset M$ to be any reducible finite index
depth 2 inclusion of von Neumann algebras with finite dimensional centers. 
A Galois correspondence has been established in the case of finite index,
finite depth inclusions of II$_1$ factors by Nikshych and Vainerman \cite{NV2,
NV3}.  The infinite index, depth 2 case has been treated by Enock and Vallin in
\cite{EV, E2} for arbitrary von Neumann  algebras endowed with a regular
operator valued weight.

\section{Quantum groupoids}

\subsection{Bialgebroids}

We try to formulate the minimal requirements on an algebraic structure 
which is to describe "symmetries", hence generalizing the notions of 
grouprings, groupoid duals, Hopf algebras, \dots etc. It must be a ring $A$ 
together with a monoidal structure on the category $\M_A$ of its right modules.
If $\M_A$ would be a bimodule category $_R\M_R$ for some ring $R$ then a
monoidal structure would be given. This motivates the

\begin{defi}  \label{bialgebroid}
Let $R^{op}\rarr{t}A\larr{s}R$ be a diagram in the category of rings\footnote{
Rings and their morphisms are always assumed unital}
such that the left and right actions of $R$ defined by $r\cdot a:=at(r)$,
$a\cdot r:=as(r)$ make $A$ into an $R$-$R$-bimodule. Equivalently, one
requires that the images of $s$ and $t$ commute in $A$. 
Then the ring and bimodule $A$ together with a monoidal structure on the
category $\M_A$ is called a bialgebroid over $R$ if the forgetful functor 
$\phi_R\colon\M_A \to \ _R\M_R $ is strictly monoidal. 
\end{defi}
This implicit but natural requirement on $\phi_R$ is equivalent to a comonoid
structure $\bra A,\gamma,\pi\ket$ in $_R\M_R$ which is compatible with the
ring structure. More precisely  
$\gamma\colon A\to A\o_R A$ and $\pi\colon A\to R$ are arrows in $_R\M_R$
such that
\bea
&&(\gamma\o_R\id)\circ\gamma\ =\ (\id\o_R\gamma)\circ\gamma\\
&&\lambda\circ(\pi\o_R\id)\circ\gamma=\id=\rho\circ(\id\o_R\pi)\circ\gamma
              \label{pi counit}\\
&&(s(r)\o_R 1)\gamma(a)\ =\ (1\o_R t(r))\gamma(a)\,,\qquad r\in R,\ a\in A
\label{pass}\\
&&\gamma(a)\gamma(b)\ =\ \gamma(ab)\,,\qquad a,b\in A\label{multi}\\
&&\gamma(1)\ =\ 1\o_R 1\,,\quad \pi(1)=1_R\\
&&\pi(t(\pi(a))b)\ =\ \pi(ab)\ =\ \pi(s(\pi(a))b)\,,\quad a,b\in A
\label{anchor}\ . 
\eea
These are essentially the same axioms as Lu's in\cite{Lu} except the different
formulation of (\ref{pass},\ref{multi},\ref{anchor}). Unfortunately quite some
explanations are needed to elucidate the meaning of Eqns (\ref{pass}) and
(\ref{multi}). The problem is that $A\o_R A$ is not a ring. It has a
sub-bimodule, however, which is. At first notice that the ring $A$ operates on
the bimodule $_RA_R$ by left multiplication so it is meaningful to write
$(a\o_R b)(a'\o_R b')$ for the result of the tensor product of intertwiners
$a\o_R b$ acting on the element $a'\o_R b'$ of the bimodule $A\o_R A$. This
convention is used in (\ref{pass}) and in the definition
$$
\Gamma:=\{\,x\in A\o_R A\,|\,(s(r)\o_R 1)x=(1\o_R t(r))x, \ r\in R\,\}
$$
which is an $R$-$R$-bimodule and a ring, too. Now (\ref{pass}) and
(\ref{multi}) just say that $\gamma\colon A\to\Gamma$ is a ring
homomorphism. 

For the equivalence of properties (1)---(6) and monoidality of the functor
$\phi_R$ it will suffice here to recall that $\phi_R$ is isomorphic to the
hom-functor $\Hom(A,\under)$ and monoidality of a hom-functor is equivalent to
a comonoid structure on the object $A$. This is how $\gamma$ and $\pi$ are
constructed.    

Definition \ref{bialgebroid} is ready made for a kind of Tannaka-Krein 
theorem which, in its weakest form, can be formulated as the

\begin{lem} \label{TK bia}
Let $\C$ be an additive monoidal category equivalent to a module category
$\M_A$ for some ring $A$. Let furthermore $F\colon\C\to\,_R\M_R$ be a strongly
monoidal\,\footnote{We use the terminology of MacLane's
\cite{MacLane} new edition!} functor to the category of bimodules over some
ring $R$. Then $A$ carries a bialgebroid structure over $R$ such that $\C$ and
$\M_A$ are equivalent as monoidal categories.
The same holds with $\M_A$ replaced by $\M_A^{\rm fgp}$, the
category of finitely generated projective $A$-modules. 
\end{lem}

In order to prepare the discussion of the "finite" bialgebroids 
let us compute here the endomorphism ring of the monoidal unit of $\M_A$.

\begin{lem} \label{End U}
If $A$ is a bialgebroid over $R$ then the monoidal unit $U_A$ of $\M_A$
is the additive group $R$ together with the $A$-action $r\triangleleft
a=\pi(s(r)a)$. Introducing the notation $Z^s:=s(R)\cap\Center A$ and
$Z^t:=t(R)\cap\Center A$ we have that $t^{-1}(Z^t)=s^{-1}(Z^s)$ is a
commutative subring $Z$ of $R$. The endomorphism ring $\End U$ is isomorphic
to $Z$ and consists of multiplications in $R$ with elements of $Z$.
\end{lem}
\begin{proof}
The action property of $\triangleleft$ follows from (6). 
Let $\xi\in\End U$. Then it is also an $R$-$R$ bimodule endomorphism of $R$,
hence $\xi(r)=rz$, $r\in R$, for some $z\in\Center R$. Hence
\beanon
as(z)&=&a\cdot z=a\1\cdot\pi(a\2)z=a\1\cdot\xi(\pi(a\2))=
a\1\cdot\xi(1_R\triangleleft a\2)\\
&=&a\1\cdot(\xi(1_R)\triangleleft
a\2)=a\1\cdot\pi(s(z)a\2)=a\1\cdot\pi(t(z)a\2)\\
&=&s(z)a\1\cdot\pi(a\2)\ =\ s(z)a
\eeanon
for all $a\in A$, therefore $s(z)\in Z^s$. Similarly, one can prove that
$t(z)\in Z^t$. Since $s$ and $t$ are sections of $\pi$, they are injective and
this proves the Lemma.
\end{proof}
The above Lemma makes it natural to consider the forgetful functor
$\phi_Z\colon \M_A\to\,_Z\M_Z$ instead of $\phi_R$. It has the advantage that
$A$ can be reconstructed from it as $\End\phi_Z$, which is not true for
$\phi_R$. However, $\phi_Z$ is a monoidal functor only in the relaxed sense.
The natural transformation 
$$ 
\mu_{V,W}\colon \phi_Z(V)\o_Z\phi_Z(W)\to\phi_Z(V\oo W) 
$$
is no longer an isomorphism but (suppressing the $\phi_Z$'s in the diagram)
\bea\label{exact}
&&\parbox[c]{4.5in}{
\begin{picture}(320,65)(5,-10) 
\put(0,40){$V\o_Z (U\o_Z W)$}
\put(0,0){$(V\o_Z U)\o_Z W$}
\put(20,35){\vector(0,-1){25}}
\put(23,22){$\wr$}
\put(125,22){$V\o_Z W$}
\put(75,40){\vector(3,-1){45}}
\put(75,5){\vector(3,1){45}}
\put(88,38){$\scriptstyle V\o_Z\lambda_W$}
\put(88,5){$\scriptstyle \rho_V\o_Z W$}
\put(165,25){\vector(1,0){50}}
\put(182,30){$\mu_{V,W}$}
\put(220,22){$V\oo W$}
\put(250,25){\vector(1,0){30}}
\put(290,22){$0$}
\end{picture}
}
\eea
is an exact sequence. Together with the inclusion map $\zeta\colon Z\to R$ the
monoidal functor $\bra\phi_Z,\mu,\zeta\ket$ contains all information about
the bialgebroid $A$. This is the content of the next

\begin{thm}
Bialgebroid structures on the ring $A$ are in one-to-one correspondence with
following categorical data.
\begin{enumerate}
\item[i)] On the one hand a monoidal structure $\bra\M_A,\oo,U\ket$ on the
category of right $A$-modules. Denoting by $Z$ the endomorphism ring of the
monoidal unit $U$, every hom-set, so $A=\End A_A$ too, becomes endowed with a
$Z$-$Z$-bimodule structure. Thus we have the forgetful functor
$\phi_Z\colon\M_A\to\,_Z\M_Z$. 
\item[ii)] On the other hand, a monoidal structure
$\bra\phi_Z,\mu,\zeta\ket$ on the forgetful functor such that (\ref{exact}) is
exact. 
\end{enumerate}
\end{thm}
\begin{proof}
The idea of the proof is to show that $\phi_Z$ factors through a strongly
monoidal forgetful functor $\phi_R\colon\M_A\to\,_R\M_R$. Here the ring 
$R$ is the additive group $U$ together with the multiplication
\bea \label{m of R}
m&\colon&U\o_Z U\rarr{\mu_{U,U}}U\oo U\rarr{\sim}U
\eea
and unit $\zeta\colon Z\to U$. Every object $V$ in $\M_A$ carries the
$R$-$R$-bimodule structure defined by the left and right actions
\bea
\lambda_V&\colon&U\o_Z V\rarr{\mu_{U,V}}U\oo V\rarr{\sim}V\\
\rho_V&\colon&V\o_Z U\rarr{\mu_{V,U}}V\oo U\rarr{\sim}V
\eea
Strong monoidality of the functor $\M_A\to\,_R\M_R$ follows from exactness of
(\ref{exact}).
\end{proof}

\subsection{Bialgebroids over separable base}
In contrast to $A\o_R A$ the bimodule tensor product $A\o_Z A$ is a ring.
This offers the tempting possibility to use, instead of $\gamma$, a
comultiplication of the  $A\to A\o_Z A$ type which could then be
multiplicative in the usual sense. For this purpose we need an embedding
$A\o_R A\subset A\o_Z A$. This is possible if $R$ is a separable algebra over
$Z$ \cite{DeMeyer-Ingraham}. In the next definition we use \cite{Abrams} to
formulate a Frobenius algebra structure as a special coalgebra structure on
the $Z$-algebra $R$. Usually an algebra being Frobenius is a {\em property},
but what we need here is a {\em structure}, i.e., a choice of a functional
$\psi$ possessing a pair of dual bases, reformulated as a comultiplication
$\delta$.

\begin{defi}
A bialgebroid over a separable base consists of a bialgebroid $\bra A, R,$ $t,
s,$ $\gamma,\pi\ket$ and of a separability structure $\bra R, Z,
\delta,\psi\ket$. The latter means that the $Z$-algebra $R$ has also a 
$Z$-coalgebra structure with $\delta\colon R\to R\o_Z R$ and $\psi\colon 
R\to Z$ which is compatible with the $Z$-algebra structure in the sense that
$\delta$ is an $R$-$R$-bimodule map 
\beq \label{Frob str}
(R\o_Z m)\circ(\delta\o_Z R)=\delta\circ m=
(m\o_Z R)\circ(R\o_Z \delta)
\eeq
(this means a Frobenius algebra structure)
and moreover 
\beq \label{sepa}
m\circ \delta=\id_R
\eeq
(separability) where $m\colon R\o_Z R\to R$ is multiplication on $R$.
\end{defi}

Notice that the data $\bra R_Z,m,\zeta,\delta,\psi\ket$ is the same as of
a bialgebra in the category of $Z$-modules, however, the compatibility
condition between the algebra and coalgebra structures is different. The
compatibility (\ref{Frob str}) does not need any symmetry or
braiding\footnote{Notice that Frobenius/separability structures could have been
defined in $_Z\M_Z$ which does not have any braiding. In fact this is what we
do but $R\in\,_Z\M_Z$ is a diagonal bimodule.} in $\M_Z$.

Once having a separability structure on $R\supset Z$ we can introduce
a natural transformation $\delta_{\scriptscriptstyle V,W}\colon V\o_R W\to
V\o_Z W$ for $R$-$R$-bimodules $V$ and $W$ that splits the canonical
epimorphism $V\o_Z W\to V\o_R W$. Namely
\bea 
\delta_{\scriptscriptstyle V,W}(v\o_R w)&:=& \sum_i \,v\cdot e_i\o_Z f_i\cdot w
\eea
where $\sum_i e_i\o_Z f_i=\delta(1_R)$. With this we can define the new
comultiplication and counit on the bialgebroid $A$. They are the
$Z$-$Z$-bimodule maps 
\bea
\cop&:=\delta_{\scriptscriptstyle A,A}\circ\gamma&\colon A\to A\o_Z A\\
\eps&:=\psi\circ\pi&\colon A\to Z\ .
\eea
These maps no longer preserve the unit, e.g. $\cop(1)=\sum_i s(e_i)\o_Z
t(f_i)$, but $\cop$ is multiplicative. The next Proposition shows how the whole
bialgebroid structure of $A$ can be reformulated in terms of $\cop$ and $\eps$
forgetting about $R$ altogether.  
\begin{pro}
A bialgebroid over separable base is equivalent to the data $\bra A, Z, t,
s,$ $\cop, \eps\ket$ where $A$ is a ring, $Z$ is a commutative ring,
$s,t\colon Z\to\Center A$ are unital ring homomorphisms making $A$ into a
$Z$-$Z$-bimodule, and $\bra \,_ZA_Z,\cop,\eps\ket$ is a comonoid in the
category $_Z\M_Z$. These data are subject to the axioms
\bea
\cop(ab)&=&\cop(a)\cop(b)\quad a,b\in A \label{Z-bia 1}\\
(\cop(1)\o_Z 1)(1\o_Z\cop(1)&=&\cop^2(1)\ =\ 
(1\o_Z\cop(1))(\cop(1)\o_Z 1) \label{Z-bia 2}\\
\eps(ab\1)\eps(b\2 c)&=&\eps(abc)\ =\ 
\eps(ab\2)\eps(b\1 c)\quad a,b,c\in A \label{Z-bia 3}\ ,
\eea
where $\cop^2$ stands for $(\cop\o_Z\id)\circ\cop\equiv(\id\o_Z\id)\circ\cop$.
Moreover, the ring $Z$ is maximal in the sense that 
\beq \label{Z-bia 4}
t(Z)=A^L\cap\Center A\,,\qquad s(Z)=A^R\cap\Center A
\eeq
where the $Z$-subalgebras $A^L$ of $_ZA$ and $A^R$ of $A_Z$ are defined by
\bea
A^L&=&\{a\in A\,|\,\cop(a)=(a\o_Z 1)\cop(1)=\cop(1)(a\o_Z 1)\,\}\\
A^R&=&\{a\in A\,|\,\cop(a)=(1\o_Z a)\cop(1)=\cop(1)(1\o_Z a)\,\}
\eea 
\end{pro}

The forgetful functor $\phi_Z\colon \M_A\to\,_Z\M_Z$ in the separable case
is not only exactly monoidal in the sense of (\ref{exact}) but is also split
in the sense of 
\begin{defi} 
The data $\bra  F, \mu,\zeta, \delta,\psi\ket$ is called
a split monoidal functor\footnote{Another possible name: "bimonoidal
functor".} if $\bra  F,\mu,\zeta\ket$ is a monoidal functor from a
monoidal category $\bra\C,\oo,u\ket$ to another $\bra\D,\o,Z\ket$, 
$\delta_{a,b}\colon  F(a\oo b)\to F(a)\o F(b)$ is a natural transformation and
$\psi$ is an arrow $ F(u)\to Z$ such that 
\begin{enumerate}
\item $\delta$ splits $\mu$, i.e., for all objects $a,b$ of $\C$ 
\beq \label{split split}
\mu_{a,b}\circ\delta_{a,b}= F(a\oo b)
\eeq
\item $\delta$ is coassociative, i.e., for all objects $a,b,c$ of $\C$ 
\beq \label{split coass}
( F(a)\o \delta_{b,c})\circ\delta_{a,b\oo c}\ =\ (\delta_{a,b}\o F(c))
\circ\delta_{a\oo b,c}
\eeq
\item $\psi$ is the counit for $\delta$, i.e.,
\beq
(\psi\o F(c))\circ\delta_{u,c}\ =\ c\ =\ 
( F(c)\o \psi)\circ\delta_{c,u}
\eeq
\item $\delta$ is compatible with $\mu$ in the sense of the equations
\bea
(\mu_{a,b}\o  F(c))\circ( F(a)\o\delta_{b.c})&=&
\delta_{a\oo b,c}\circ\mu_{a,b\oo c} \label{split comp 1}\\
( F(a)\o\mu_{b,c})\circ(\delta_{a,b}\o F(c))&=&
\delta_{a,b\oo c}\circ\mu_{a\oo b,c} \label{split comp 2}
\eea
\end{enumerate}
\end{defi}

Notice that equations (\ref{split coass}, \ref{split comp 1}, \ref{split comp
2}) are just variations of the associativity condition on $\mu$ in which
certain $\mu$ arrows were replaced with oppositely oriented $\delta$'s.  
These equations have interesting similarity with the axioms of a Frobenius
structure (\ref{Frob str}) while the splitting property (\ref{split split})
corresponds to the separability axiom (\ref{sepa}). 

Split monoidal functors are just the functors arising as forgetful
functors $\M_A\to\,_Z\M_Z$ for a bialgebroid over separable base. We give here
the precise statement for representable functors.

\begin{thm}
Let $\bra \C,\oo, u\ket$ be a monoidal category with finite progenerator $g$.
With the notation $Z=\End u$ let $ F\colon\C\to\,_Z\M_Z$ be the hom-functor
$ F=\Hom(g,\under)$. Then split monoidal structures
$\bra F,\mu,\zeta,\delta,\psi\ket$ for $ F$ are in bijective
correspondence with bialgebroid structures $\bra _RA_R,\gamma,\pi\ket$ on the
ring $A=\End g$ that have separable base. Moreover, any such split monoidal
$ F$ factorizes, as a monoidal functor, through an equivalence
$\C\cong\M_A^{fgp}$ of monoidal categories to the categegory of finitely
generated projective right $A$-modules.
\end{thm}
\begin{proof}
Since $g$ is a finite progenerator, $\C\cong\M_A^{fgp}$. Monoidality of $ F$
determines the comonoid $\bra g,\gamma,\pi\ket$ where
\bea
\gamma\colon g\to g\oo g\,,&&\quad \gamma:=\mu_{g,g}(g\o_Z g)\\
\pi\colon g\to u\,,&&\quad\pi:=\zeta(u)\ .
\eea
Using also the split monoidality structure we can define
\bea
\cop\colon A\to A\o_Z A\,,&&\quad \cop(a):=\delta_{g,g}(\gamma\circ a)\\
\eps\colon A\to Z\,,&&\quad\eps(a):=\psi(\pi\circ a)
\eea
which form a comonoid $\bra A,\cop,\eps\ket$ in $_Z\M_Z$ and can be verified to
obey the properties (\ref{Z-bia 1}, \ref{Z-bia 2}, \ref{Z-bia 3}). The
maximality property (\ref{Z-bia 4}) holds automatically by the very definition
of $Z$ as $\End u$. Thus $A$ is a bialgebroid with separable base $R=\Hom(g,u)$
endowed with multiplication as in (\ref{m of R}). 
\end{proof}

Antipodes can be introduced on bialgebroids by postulating the existence of
left and right dual objects in the category $\M_A$. This will not be discussed
here.

\subsection{Weak Hopf algebras}

The (\ref{Z-bia 1}, \ref{Z-bia 2}, \ref{Z-bia 3}) axioms are already very 
close to the weak bialgebra axioms of
\cite{BNSz}. In fact a weak bialgebra (WBA) arises from a bialgebroid over
separable base by a further finiteness condition: The commutative ring $Z$
should be a separable algebra over a field $K$. Then the structure maps can be
formulated in the symmetric monoidal category of $K$-vector spaces. Denoting
the tensor product over $K$ by $\o$ there is a comultiplication $\cop\colon
A\to A\o A$ and counit $\eps\colon A\to K$ satisfying exactly the axioms
(\ref{Z-bia 1}, \ref{Z-bia 2}, \ref{Z-bia 3}) except that $Z$ is replaced 
everywhere with $K$. Then $Z$, more precisely two copies of it, $s(Z)=Z^R$
and $t(Z)=Z^L$, can be reconstructed as $Z^c=A^c\cap\Center A$, $c=L,R$,
respectively. So the (\ref{Z-bia 4}) maximality condition 
is not needed, although it might be reasonable to demand that the ground
field be intrinsically defined by the bialgebroid structure. This could be
achieved by adding the axiom that the hypercenter $Z^L\cap Z^R$ of $A$ is a
field $K$. (This kind of weak bialgebras (and weak Hopf algebras) are called
indecomposable.)  

A weak Hopf algebra (WHA) over $K$ is a WBA $A$ over $K$ such that there 
exists a linear map $S\colon A\to A$, called the antipode, such that 
\bea
a\1 S(a\2)&=&\PL(a)\\
S(a\1)a\2&=&\PR(a)\\
S(a\1)a\2S(a\3)&=&S(a)
\eea
for all $a\in A$ where $\PL$, $\PR$ are analogues of $\pi$ in (\ref{pi
counit}) and are defined by $\PL(a)=\eps(1\1 a)1\2$, $\PR(a)=1\1\eps(a1\2)$.
The antipode, if exists, is unique. It is antimultiplicative,
anticomultiplicative and maps $A^L$ onto $A^R$ bijectively. In the sequel we
shall also assume that $A$ is finite dimensional over $K$. In this case $S$ is
invertible. 

The dual space $\duA=\Hom_K(A,K)$ of a WHA endowed with multiplication
and comultiplication obtained by transposing the comultiplication and
multiplication of $A$, respectively, is again a WHA over $K$. Moreover
there is a natural identification of their left, right subalgebras: $A^L\iso
\duA^R$, $l\mapsto l\la\du1$ and $A^R\iso\duA^L$, $r\mapsto \du1\ra r$ given
by the Sweedler arrows: For $a\in
A$ and $\varphi\in\duA$ one writes $a\la\varphi:=\varphi\1\bra\varphi\2,a\ket$
and $\varphi\ra a:=\bra\varphi\1,a\ket\varphi\2$. 

\begin{defi} \label{def: integral} 
A left (right) integral in a weak Hopf algebra $A$ is an 
element $\iL\in A$ ($\iR\in A$) satisfying 
\beq                                                 \label{l (r)} 
x\iL\ =\ \PL(x)\iL\qquad \left(\ \iR x\ =\ \iR\PR(x)\ \right) 
\eeq 
for all $x\in A$.
$\iL$ is called normalized if $\PL(\iL)=1$ ($\iR$ is called 
normalized if $\PR(\iR)=1$). 
A left or right integral in $A$ is called non-degenerate if it 
defines a non-degenerate functional on $\duA$. 
\end{defi} 

Existence of non-zero integrals follows from a theorem on weak Hopf modules
(Thm 3.9 of \cite{BNSz}). 
The existence of non-degenerate or normalized integrals in a WHA are related
to the $K$-algebra $A$ to be Frobenius or semisimple, respectively \cite{BNSz}.
For example the next result provides a weak Hopf version of Maschke's Theorem.
  
\begin{thm}[see \cite{BNSz} Thm 3.13]       \label{thm: semi} 
The following conditions on a WHA $A$ over $K$ are equivalent: 
\begin{enumerate} 
\item[i)] $A$ is semisimple. 
\item[ii)] There exists a normalized left integral $\iL\in A$. 
\item[iii)] $A$ is a separable $K$-algebra. 
\end{enumerate} 
\end{thm} 

\begin{defi}  \label{def: Haar} 
An element $h$ of a WHA $A$ is called a Haar integral in $A$  
if $h$ is a normalized 2-sided integral, 
i.e., $h$ is a left integral, a right integral, and $\PL(h)=\PR(h)=1$. 
\end{defi} 

\begin{thm}[see \cite{BNSz} Thm. 3.27]         \label{thm: Haar} 
Let $A$ be a WHA over an algebraically closed field $K$. Then a 
necessary and sufficient condition for the existence of Haar 
integral $h\in A$ is that $A$ is semisimple and there exists an invertible
element $g\in A$ such that $gxg^{-1}=S^2(x)$ for $x\in A$ 
and $\tr_r(g^{-1})\neq 0$ in all irreducible representations 
$r$ of $A$\,. 
\end{thm} 
If exists, the Haar integral is unique and is an idempotent.

\section{$C^*$-weak Hopf algebras}

A $C^*$-weak Hopf algebra is a WHA $A$ over $\mathbb{C}$ which is a finite
dimensional $C^*$-algebra and the comultiplication $\cop$ is a $^*$-algebra
map. By uniqueness of the antipode it follows that $S(S(a)^*)^*=a$ for all
$a\in A$. If also $S^2=\id$, the $C^*$-WHA is called a weak Kac algebra
\cite{NV1}.  The counit $\eps\colon A\to \mathbb{C}$ is always a positive
linear functional and the associated GNS representation is the monoidal unit
$U$ of the representation category $\M_A$. If $U$ is irreducible,
or equivalently, if $Z^L=\mathbb{C}1$, i.e., the inclusion $A^L\subset A$ is
connected, then $A$ is called pure \cite{BSz} or connected \cite{NV1}. 

\subsection{The Haar measure} 
\begin{thm}             \label{thm: C^* Haar} 
In a $C^*$-WHA $A$ there exists a Haar integral. It is 
a selfadjoint $S$-invariant idempotent, $h=h^*=h^2=S(h)\in A$, such that 
\beq 
(\varphi,\psi)\ :=\ \bra\varphi^*\psi,\ h\ket\ , 
\qquad \varphi,\psi\in\duA\ , 
\eeq 
is a scalar product making $\duA$ a Hilbert space and 
making the left regular module $_{\duA}\duA$ a faithful 
$^*$-representation of the $^*$-WHA $\duA$. Thus $\duA$ is a 
$C^*$-WHA, too. 
\end{thm} 

Thus also $\duA$ has a Haar integral $\hat h\in \duA$.
This provides the faithful conditional expectations 
\bea   \label{Haar cond} 
E^L\colon A\to A^L\,,&\quad&E^L(x)=\hat h\la x\\ 
E^R\colon A\to A^R\,,&\quad&E^R(x)=x\ra \hat h 
\eea 
It can be shown that $\hat h\la h\in A^L$ and $h\ra\hat h\in A^R$ are positive
and invertible. The so called canonical grouplike element is defined by
\bea 
g&:=& g_Lg_R^{-1}\ ,\label{eq: g=}\\ 
g_L&:=&(\hat h\la h)^{1/2}\ ,\quad g_R\ =\ (h\ra\hat h)^{1/2} 
\eea 
and can be characterized as the unique $g\in A$ such that 
\begin{enumerate} 
\item[i)] $g\geq 0$ and invertible, 
\item[ii)] $gxg^{-1}=S^2(x)$ for all $x\in A$, 
\item[iii)] $\tr_r(g^{-1})=\tr_r g$ in all irreducible representations $r$.
\end{enumerate} 
In general the Haar functional $\bra \hat h,\under\ket\colon A\to\mathbb{C}$ is
not a trace but instead
\beq
\bra\hat h,ab\ket\ =\ \bra\hat h,b\,g_Lg_Ra(g_Lg_R)^{-1}\ket\quad a,b\in A\,.
\eeq
It is a trace iff $S^2=\id$, i.e., iff $A$ is a weak Kac algebra.

\subsection{Dimensions}
The category $\rep A$ of finite dimensional $^*$-representations of a
$C^*$-WHA $A$ is a monoidal category with monoidal structure inherited from
the forgetful functor to the category of Hilbert $A^L$-$A^L$-bimodules, a
$^*$-functor analogue of the $\phi_R$ of Section 1. Since the usual convention
in $^*$-representations is left action, the functor is constructed by
considering $A$ to be a bimodule via $l_1\cdot a\cdot l_2:=l_1 S^{-1}(l_2)a$,
$l_1,l_2\in A^L$, $a\in A$. Then monoidal product of two representations
$D_i\colon A\to\B(H_i)$, $i=1,2$, is
defined on the $A^L$-$A^L$-Hilbert bimodule tensor product 
$H_1\o_{A^L}H_2$ endowed with the
left action via the comultiplication, $D_1\oo D_2:=(D_1\o D_2)\circ\cop$, which
is well defined due to the identities (2.31a-b \cite{BNSz}). The monoidal unit
of $\rep A$ is the GNS representation $D_\eps$ associated to the counit
$\eps\colon A\to\mathbb{C}$. $D_\eps$ is irreducible iff $A$ is pure.

All objects $D$ of $\rep A$ have conjugates $\bar D$, i.e., two-sided duals,
defined by help of the antipode \cite{BSz2}. If the WHA is pure then all
conditions are fulfilled to apply the theory of dimensions of \cite{LR}. Even
if $A$ is not pure one can find analogues of the {\em standard}
conjugacy intertwiners $R_D\colon D_\eps\to\bar D\oo D$, $\bar R\colon D_\eps
\to D\oo \bar D$. If $D$ is irreducible then $R^*_D\circ R_D$ and $\bar
R^*_D\circ\bar R_D$ are selfintertwiners of $D_\eps$ proportional to one
minimal projection in $\End D_\eps=D_\eps(Z^L)$ (cf. Lemma \ref{End U}). But
not to the same projection, in general. Standard normalization means 
choosing $R_D$, $\bar R_D$ for all objects so that it respects direct sums,
like in \cite{LR}, and for irreducible objects $D$
\beq
R^*_D\circ R_D=d_D D_\eps(z^L_\mu)\,,\quad
\bar R^*_D\circ\bar R_D=d_D D_\eps(z^L_\nu)
\eeq
with the same positive (in fact $\geq 1$) number $d_D$ in both equations, but
with possibly different minimal projections $z^L_\mu,z^L_\nu\in Z^L$. All
these data on the right hand sides depend only on the equivalence class $q$ to
which $D$ belongs. The number $d_q=d_D$ is called the dimension of the sector
$q$, while $q^L=\nu$ and $q^R=\mu$ are called the left and right vacuum of $q$,
respectively. The dimensions of irreducibles can be written as
\beq
d_q= k(q^L)^{-1/2}k(q^R)^{-1/2}\ \tr_q\, g\,,\quad k(\mu):=\eps(z^L_\mu)\,.
\eeq

For pure WHA's there is only one vacuum sector. This is the case when
$D\mapsto d_D$ is an additive and multiplicative dimension function. For
general $C^*$-WHA's one forms the matrix $\d_q=d_qe_{\nu\mu}$ (a number
times a matrix unit) for all sectors $q$, the rows and columns of which are
labelled by the set of vacua, i.e., by the irreducibles contained in $D_\eps$.
For an arbitrary representation $D$ one defines the matrix $\d_D:=\sum_q
N_q(D)\d_q$, where $N_q(D)$ is the multiplicity of $q$ in $D$. The so defined
dimension matrix will then be both additive and multiplicative. Conjugating the
representation its dimension matrix goes to its transposed matrix. 

Particularly interesting is the dimension matrix $\d_A$ of the left regular
representation. It turns out to be similar to the matrix $\d_{\duA}$, which
is computed, of course, in another category, in $\rep \duA$. But there exists
a matrix $\d^L$ with non-negative coefficients, and its transposed matrix
$\d^R$, such that $\d_A=\d^L\d^R$ and $\d_{\duA}=\d^R\d^L$. These new matrices
can be interpreted as the dimension matrices of $A^L$ and $A^R$,
respectively \cite{BSz2}. 

In the next theorem we assume that
$A$ is an indecomposable $C^*$-WHA, i.e., $Z^L\cap Z^R=\mathbb{C}1$. 

\begin{thm}[see \cite{BSz2}] \label{Markov}
The basic constructions for the inclusions $A^L\subset A$ and
$\duA\supset\duA^R$ coincide and equal to the smash product $C^*$-algebra
$A\# \duA$.

There exists a unique normalized trace $\tau$, called the Markov trace, on the
smash product such that for all vacuum $\nu$ of $A$ and all vacuum $\hat\nu$
of $\duA$ the restrictions
$$
\tau\upharpoonright z^L_\nu A\,,\quad
\tau\upharpoonright z^R_\nu A\,,\quad
\tau\upharpoonright \hat z^L_{\hat\nu}\duA\,,\quad
\tau\upharpoonright \hat z^R_{\hat\nu}\duA
$$
are the Markov traces of the connected inclusions 
$$
z^L_\nu A^L\subset z^L_\nu A\,,\quad
z^R_\nu A^R\subset z^R_\nu A\,,\quad
\hat z^L_{\hat\nu}\duA^L\subset \hat z^L_{\hat\nu}\duA\,,\quad
\hat z^R_{\hat\nu}\duA^R\subset \hat z^R_{\hat\nu}\duA\,,
$$
respectively. 
The corresponding trace preserving conditional expectations all have the same
index $\delta$. This index coincides with the common Perron-Frobenius
eigenvalue of the regular dimension matrices $\d_A$ and $\d_{\duA}$. 
\end{thm}

If $A\cong\oplus_q M_{n_q}(\mathbb{C})$ is pure one
gets the number $\delta=\d_A=\sum_qn_qd_q$.   

The Markov conditional expectations $A\to A^L$, $A\to A^R$ are different from
the Haar conditional expectations (\ref{Haar cond}), unless $A$ is a weak Kac
algebra. In this latter case $\delta$ is an integer. The Haar conditional
expectations $E^L$ and $E^R$ also have a common scalar index $I$, but $I\geq
\delta$, in general. 

\begin{exa}
In \cite{BSz} we gave an example of a $C^*$-WHA structure on the matrix
algebra $A=M_2\oplus M_3$. The two sectors obey the fusion rules $3\x 3=2+3$,
with $2$ being the unit of the fusion ring. $A^L\cong M_1\oplus M_1$ and
$A^L\cap\,\Center A=A^L\cap A^R=\mathbb{C}1$, so both $A$ and $\duA$ are
connected (i.e., $A$ is biconnected). What is more $A\cong\duA$. The
dimensions of the sectors are $d_2=1$, $d_3=(1+\sqrt{5})/2$. The Haar index is
$I=4+2d_3=5+\sqrt{5}=7.24$ and the Markov index is $\delta=2+3d_3=6.85$.
\end{exa}

The above example is the first of a series of WHA's with the underlying
algebra being a Temperley-Lieb algebra \cite{NV3, NV4}.

Finally we mention that there exists a description of weak $C^*$-Hopf algebras
in terms of finite dimensional multiplicative partial isometries \cite{BSz3,
Val2}.

\section{Finite index depth 2 inclusions}
The most important (sofar the only) application of $C^*$-WHA's is
the characterization of finite index, depth 2 inclusion of von Neumann
algebras. 

\subsection{Weak Hopf actions}
A left action of a $C^*$-WHA $A$ on the unital $C^*$-algebra $M$ is an algebra
map $\alpha\colon A\to\End_\mathbb{C}M$ which respects the $^*$-algebra
structure,
\bea
\alpha_a(mm')&=&\alpha_{a\1}(m)\alpha_{a\2}(m')\,,\quad a\in A,\ m,m'\in M\\
\alpha_a(m)^*&=&\alpha_{S(a)^*}(m^*)
\eea
and leaves the identity "invariant" in the sense of the relation
\beq \label{action on 1}
\alpha_a(1_M)\ =\ \alpha_{\pi^L(a)}(1_M)\ .
\eeq
One also requires that $m\mapsto \alpha_a(m)$ is continuous for all $a\in A$.
The invariants of a left action are the elements $n\in M$ which transform like
the identity in (\ref{action on 1}). The invariants form a $C^*$-subalgebra
$M^A$ which can be expressed as the result of the application of the Haar
integral, $M^A=\alpha_h(M)$. 

Since the trivial representation of a WHA is not one-dimensional, together with
$1_M$ one should consider on equal footing all operators $\alpha_l(1_M)$,
$l\in A^L$. These operators form a $^*$-subalgebra $M^R$ and $l\mapsto
\alpha_l(1_M)$ is a $^*$-algebra epimorphism. (For faithful actions it is an
isomorphism.) One considers $M$ as a right $A^L$-module by setting $m\cdot
l=m\alpha_l(1_M)$.

The crossed product $C^*$-algebra $M\rtimes A$ can be defined as the
universal $C^*$-algebra of the $^*$-algebra defined on the bimodule tensor
product $M\o_{A^L}A$ by the relations
\bea
(m\rtimes a)(n\rtimes b)&=&m\alpha_{a\1}(n)\rtimes a\2b\\
(m\rtimes a)^*&=&\alpha_{a\1^*}(m^*)\rtimes a\2^*
\eea
The subalgebras $1_M\rtimes A$ and $M\rtimes 1$ will be identified with $A$
and $M$, respectively. One is interested in situations when the triple
$M^A\subset M\subset M\rtimes A$ is a basic construction. For that we have to
select a class of "nice" actions.
\begin{defi}[\cite{NSzW}]
A regular action of $A$ on $M$ is an action $\alpha$ such that
\begin{enumerate}
\item[i)] $M^R = A^L$,
\item[ii)] $M'\cap(M\rtimes A) = A^R$, and
\item[iii)] $\alpha_h\colon M\to M^A$ is a conditional expectation of finite
index \cite{Watatani}.
\end{enumerate}
\end{defi}
Regular actions are outer in the sense of the relative commutant
$M'\cap(M\rtimes A)$ being as small as possible. ($A^R$ always commutes with
$M$ in the crossed product.) For the more fundamental meaning of outerness, as
opposite to (partly) inner actions, see \cite{NSzW}. 
\begin{pro}[\cite{NSzW}]
If $\alpha$ is a regular action of $A$ on the $C^*$-algebra $M$ then
\begin{enumerate}
\item The crossed product $M_2:=M\rtimes A$ is the basic construction for the
finite index inclusion $N:=M^A\subset M$.
\item $N'\cap M=A^L$
\item $M'\cap M_2=A^R$
\item $N'\cap M_2=A$
\item $\Center N=Z^L$, $\Center M=A^L\cap A^R$, $\Center M_2=Z^R$.
\end{enumerate}
\end{pro}
Regular actions are Galois actions: Denoting by $\rho\colon M\to M\o_{A^L}\duA$
the right coaction associated to the left action $\alpha$, the canonical map
\beq
M\o_N M\to M\o_{A^L}\duA\,,\quad (m\o m')\mapsto (m\o \hat 1)\rho(m')
\eeq
is an isomorphism. For the proof see the Appendix of \cite{NSzW}.

\subsection{Depth 2 inclusions of II$_1$ factors}
This subsection is based on the results obtained by D. Nikshych and L.
Vainerman in \cite{NV2, NV3}. Let $N\subset M$ be a finite index, depth 2 
inclusion of II$_1$ factors and let $E\colon M\to N$ denote the trace
preserving conditional expectation. Then one constructs the Jones tower
$$
N\stackrel{E_1}{\subset}M\stackrel{E_2}{\subset} M_2\stackrel{E_3}{\subset}
M_3\subset\dots
$$
with the finite index conditional expectations $E_n\colon M_n\to M_{n-1}$
implemented by Jones projections $e_n\in M'_{n-1}\cap M_{n+1}$ satisfying the
Temperley-Lieb algebra with $e_n e_{n+1} e_n= e_n/\delta$, where $\delta$ is
the index of $E$, i.e., the minimal index of $N\subset M$. The derived tower
$\dots \subset N'\cap M_n\subset N'\cap M_{n+1}\subset\dots$ consists of
finite dimensional $C^*$-algebras and is a Jones tower starting from $N'\cap
M$ (depth 2 condition). 

Define the algebras $A:=N'\cap M_2$, $B:=M'_1\cap M_3$ and a pairing, i.e.,
a non-degenerate bilinear form (cf Eqn 4.134 in \cite{BSz2})
\beq \label{pair}
\bra a,b\ket:=\delta^{3/2}\tau(e_2e_1 b \cL a \cR)\,,
\quad a\in A,\ b\in B\,,
\eeq
where $\cL\in \Center(N'\cap M)$, $\cR\in \Center(M'\cap M_2)$ are positive
invertible elements.
The pairing yields the coalgebras $\bra A,\cop_A,\eps_A\ket$ and $\bra
B,\cop_B,\eps_B\ket$ immediately. Antipodes can be introduced by 
$$ S_A\colon A\to A\,,\ 
S_A(a):=(a_*)^*\,\quad\mbox{where}\ \bra a_*,b\ket:=\overline{\bra a, b^*\ket}
$$
and an analogue expression for $S_B\colon B\to B$.
The difficult part is to show that for an appropriate (in fact unique) choice
of $\cL$, $\cR$ the comultiplication $\cop_B$ is a $^*$-algebra map. If it is
done all other axioms of weak Hopf algebras hold automatically. 
The next Theorem is a reformulation of the results of D. Nikshych and 
L. Vainerman and provides a 
generalization of the duality theorem of irreducible inclusions of
factors and $C^*$-Hopf algebra actions \cite{Szymanski, Longo, David}.  
\begin{thm}[cf. \cite{NV2}]
Let $N\subset M\subset M_2\subset M_3$ be the Jones tower built on a finite
index, depth2 inclusion $N\subset M$ of II$_1$ factors. Then there are
biconnected $C^*$-WHA structures on the relative commutants $A=N'\cap M_2$ and
$B=M'\cap M_3$ such that 
\begin{enumerate}
\item $A$ and $B$ are the duals of each other w.r.t. the pairing (\ref{pair})
\item $A$ acts regualarly on $M$ with $N$ being the invariant subalgebra
\item $M_2$ is isomorphic to the crossed product $M\rtimes A$
\item the index $[M:N]$ is equal to the Markov index $\delta$ of the
finite dimensional incusion $A^L\subset A$ determined by Thm \ref{Markov}. 
\end{enumerate}
The WHA's become unique under the additional requirement of $S^2\upharpoonright
A^L$ be the identity. 
\end{thm} 
Although the generalization to non-factors is quite plausible and has already
been suggested in \cite{NSzW}, no published results are available yet. There
is, however, another approach by M. Enock and J.-M. Vallin by means of which
almost arbitrary depth 2 inclusions can be described as invariant subalgebras
w.r.t. actions of Hopf bimodules \cite{EV, E2}. 

In \cite{NV3} Nikshych and Vainerman made an important step in another
direction. They considered any finite depth inclusion $N\subset M$ of II$_1$
factors. Since there is always a depth 2 subtower $N\subset M_p\subset M_{2p}
\subset\dots$ of the Jones tower over $N\subset M_1$, there is always a WHA $A$
acting on $M_p$ with $N$ being its invariant subalgebra. Let its dual weak
Hopf algebra $M_p'\cap M_{3p}$ be denoted by $B$. The question arises how the
intermediate factors $M_k$, $2p< k < 3p$ are related to substructures of
$B$?  The appropriate substructure is a {\em left coideal $^*$-subalgebra},
i.e., a unital $^*$-subalgebra $I\subset B$ such that $\cop(I)\subset B\o I$.
These coideal subalgebras form a lattice with minimal element $B^L$ and maximal
element $B$. The next theorem extends the earlier results of \cite{ILP, E3}. 
\begin{thm}[see \cite{NV3}]
Let $N\subset M_1\subset M_2\subset M_3\subset\dots$ be the tower constructed
from a depth 2 subfactor $N\subset M_1$ and let $B=M_1'\cap M_3$ be the
corresponding quantum groupoid. Then intermediate factors $M_2\subset K\subset
M_3$ and left coideal $^*$-subalgebras $I\subset B$ are in one-to-one
correspondence via 
$$
K\mapsto I=M'_1\cap K\subset B\quad\mbox{and}\quad I\mapsto K=M_2\rtimes
I\subset M_3\ .
$$
\end{thm}

\subsection{Abstract inclusions}
The above results are very probably only special cases of a much more general
duality between inclusions and quantum groupoid actions. I try to outline here
the construction of bialgebroids from depth 2 arrows in a 2-category. The data
what one needs for this construction to work is reminiscent to the data of an
abstract Q-system proposed by Longo in \cite{Longo}, although we work in the
non-$^*$ framework.

Let $a\colon M\to N$ be an arrow in an additive 2-category $\C$. Assume $a$
has a left dual $a^L\colon N\to M$ with unit $\eta\colon N\to a\oo a^L$ and
counit $\eps\colon a^L\oo a\to M$. (One may think $\C$ to be the 2-category of
categories and $a$ to be the forgetful functor corresponding to a ring
inclusion $N\to M$. Then $a^L$ is the induction functor.)

Then $g:=a^L\oo a$ has a comonoid structure in the monoidal category $\C(M,M)$
defined by
\bea
&\delta:=a^L\oo \eta\oo a&\colon g\to g\oo g\\
&\pi:=\eps&\colon g\to M
\eea
Let $F$ denote the hom-functor $\Hom(g,\under)$ from $\C(M,M)$ to the category 
$_R\M_R$ of $R$-$R$-bimodules where $R:=\End a$ and the bimodule structure on
a hom-group is defined by
\beq
r\cdot h:=h\circ(r^L\oo a)\,,\qquad h\cdot r:=h\circ(a^L\oo r)\,.
\eeq
Here $r\mapsto r^L$ denotes the action of the left dual functor mapping $R$ to
$\End a^L$. Clearly, $F$ factorizes through the category of right $A:=\End
g$-modules:
\beq
\parbox[c]{4in}{
\begin{picture}(180,70)
\put(20,50){$\C(M,M)$}
\put(70,52){\vector(1,0){70}}  \put(102,55){$F$}
\put(150,50){$_R\M_R$}
\put(40,40){\vector(2,-1){47}} 
\put(92,10){$\M_A$}
\put(110,17){\vector(2,1){50}} \put(138,21){$\phi$}
\end{picture}
}
\eeq
where the forgetful functor $\phi$ is meant w.r.t. the bimodule structure 
on $A$ as in Definition \ref{bialgebroid} with source and target maps
\beq
t(r)\ :=\ r^L\oo a\,,\qquad s(r)\ :=\ a^L\oo r\,.
\eeq
In order to demonstrate that $\phi$ is monoidal we need the natural
transformation
\bea
\mu_{b,c}&\colon& F(b)\o_R F(c)\to F(b\oo c)\\
\mu_{b,c}(x\o_R y)&=&(x\oo y)\circ\delta
\eea
which is well-defined due to $(r\oo a^L)\circ\eta=(a\oo r^L)\circ\eta$, and it
is an $R$-$R$-bimodule map. Together with the arrow
\beq
\nu\colon R\to F(M)\,,\quad r\mapsto \eps\cdot r=r\cdot\eps=\eps\circ(a^L\oo r)
\eeq
$\mu$ satisfies associativity and unit constraints establishing the monoidal
functor $\bra\phi,\mu,\nu\ket$.

$\nu$ is in fact an isomorphism. In order for $\mu$ to be also an isomorphism
we have to make a further assumption on the arrow $a$ and restrict $F$ to an
appropriate subcategory of $\C(M,M)$. Assume that $a$ is of {\em depth 2},
i.e., $a\oo a^L\oo a$ is a direct summand of a finite direct sum of $a$'s. Let
$\C_g$ be the full subcategory of $\C(M,M)$ the objects $b$ of which are
direct summands of a finite multiple of $g$'s. Then $\C_g$ contains the tensor
powers $g^{\oo n}$, is a subcategory closed under the monoidal product,
and has subobjects.

We need yet a further assumption, namely that $M$ belongs to $\C_g$. This is
equivalent to the assumption that $M$ is contained in $g$ as a direct summand.

\begin{thm}
Let $\C$ be an additive 2-category\,\footnote{In fact the theorem holds for
bicategories, too.} closed under direct sums and subobjects of arrows
($=$ 1-cells) and let the arrow $a\colon N\to M$ possess a left dual, be of
depth 2, and be such that $M$ is contained in $g=a^L\oo a$ as a direct
summand. Then there is a full monoidal subcategory $\C_g$ of $\C(M,M)$ which
is equivalent,  as a monoidal category, to the category  $\M_A^{\rm f.g.p.}$ of
finitely generated projective right $A$-modules for a uniquely determined
bialgebroid structure on $A=\End g$ over the base $R=\End a$. The functor
$\Hom(g,\under)\colon\C_g\to\,_R\M_R$ factorizes, as a monoidal functor,
through the forgetful functor $\phi\colon\M_A\to\,_R\M_R$.  
\end{thm}
\begin{proof}
Need to show that $\mu_{b,c}$ is an isomorphism for objects $b$, $c$ of $\C_g$.
Choosing a direct sum decomposition $a\oo b\rarr{e_i}a\rarr{f_i}a\oo b$ we can
explicitely write down the inverse as
$$
\mu_{b,c}^{-1}(t)=\sum_i\ (\eps\oo b)\circ(a^L\oo f_i)\ \o_R\ 
                   (\eps\oo c)\circ (a^L\oo e_i\oo c)\circ
                   (g\oo t)\circ(a^L\oo\eta\oo a)
$$
for $t\in\Hom(g,b\oo c)$. This proves strong monoidality of $\bra
F\upharpoonright\C_g,\mu,\nu\ket$. From the construction of $\C_g$ it is clear
that $\C_g\cong \M_A^{\rm f.g.p.}$, as categories. Use this equivalence to
define a monoidal structure on $\M_A^{\rm f.g.p.}$. Now apply
Lemma \ref{TK bia} to conclude that $A$ is a bialgebroid over $R$ and the
equivalence is that of monoidal categories. The monoidal factorization through
$\phi$ holds by the very definition of the monoidal structure of $\M_A^{\rm
f.g.p.}$. 
\end{proof}

\end{document}